\documentclass{amsart}

\usepackage{amsmath,amsthm,amsfonts,amscd,amssymb}

\newtheorem{thm}{Theorem}[section]
\newtheorem{cor}[thm]{Corollary}
\newtheorem{lem}[thm]{Lemma}

\theoremstyle{definition}

\theoremstyle{remark}

\begin{document}

\title{Siegel's Lemma with Additional Conditions}
\author{Lenny Fukshansky}

\address{Department of Mathematics, 3368 TAMU, Texas A\&M University, College Station, Texas 77843-3368}
\email{lenny@math.tamu.edu}

\begin{abstract}
Let $K$ be a number field, and let $W$ be a subspace of $K^N$, $N \geq 1$. Let $V_1,...,V_M$ be subspaces of $K^N$ of dimension less than dimension of $W$. We prove the existence of a point of small height in $W \setminus \bigcup_{i=1}^M V_i$, providing an explicit upper bound on the height of such a point in terms of heights of $W$ and $V_1,...,V_M$. Our main tool is a counting estimate we prove for the number of points of a subspace of $K^N$ inside of an adelic cube. As corollaries to our main result we derive an explicit bound on the height of a non-vanishing point for a decomposable form and an effective subspace extension lemma.
\end{abstract}

\maketitle

\def\A{{\mathcal A}}
\def\B{{\mathcal B}}
\def\C{{\mathcal C}}
\def\D{{\mathcal D}}
\def\F{{\mathcal F}}
\def\x{{\mathcal H}}
\def\I{{\mathcal I}}
\def\J{{\mathcal J}}
\def\K{{\mathcal K}}
\def\L{{\mathcal L}}
\def\M{{\mathcal M}}
\def\R{{\mathcal R}}
\def\s{{\mathcal S}}
\def\V{{\mathcal V}}
\def\X{{\mathcal X}}
\def\Y{{\mathcal Y}}
\def\H{{\mathcal H}}
\def\cee{{\mathbb C}}
\def\pee{{\mathbb P}}
\def\que{{\mathbb Q}}
\def\real{{\mathbb R}}
\def\zed{{\mathbb Z}}
\def\aaa{{\mathbb A}}
\def\qbar{{\overline{\mathbb Q}}}
\def\eps{{\varepsilon}}
\def\ahat{{\hat \alpha}}
\def\bhat{{\hat \beta}}
\def\gt{{\tilde \gamma}}
\def\h{{\tfrac12}}
\def\be{{\boldsymbol e}}
\def\bei{{\boldsymbol e_i}}
\def\bc{{\boldsymbol c}}
\def\bm{{\boldsymbol m}}
\def\bk{{\boldsymbol k}}
\def\bi{{\boldsymbol i}}
\def\bl{{\boldsymbol l}}
\def\bq{{\boldsymbol q}}
\def\bu{{\boldsymbol u}}
\def\bt{{\boldsymbol t}}
\def\bs{{\boldsymbol s}}
\def\bv{{\boldsymbol v}}
\def\bw{{\boldsymbol w}}
\def\bx{{\boldsymbol x}}
\def\bX{{\boldsymbol X}}
\def\bz{{\boldsymbol z}}
\def\bwy{{\boldsymbol y}}
\def\bY{{\boldsymbol Y}}
\def\bL{{\boldsymbol L}}
\def\ba{{\boldsymbol\alpha}}
\def\bb{{\boldsymbol\beta}}
\def\bet{{\boldsymbol\eta}}
\def\bxi{{\boldsymbol\xi}}
\def\bo{{\boldkey 0}}
\def\bol{{\boldkey 1}_L}
\def\ep{\varepsilon}
\def\p{\boldsymbol\varphi}
\def\q{\boldsymbol\psi}
\def\rank{\operatorname{rank}}
\def\aut{\operatorname{Aut}}
\def\lcm{\operatorname{lcm}}
\def\sgn{\operatorname{sgn}}
\def\spn{\operatorname{span}}
\def\md{\operatorname{mod}}
\def\Norm{\operatorname{Norm}}
\def\dim{\operatorname{dim}}
\def\det{\operatorname{det}}
\def\Vol{\operatorname{Vol}}
\def\rk{\operatorname{rk}}

\section{Introduction and notation}

The name Siegel's Lemma is usually used to denote results about small-height solutions of a system of linear equations. Such a result in a simple form was first proved by Thue in 1909 (\cite{thue}, pp. 288-289) using the Dirichlet's box principle. Siegel (\cite{siegel}, Bd. I, p. 213, Hilfssatz) was the first to formally state this principle in the classical case. 
\smallskip

Notice that a small-height solution to a system of linear equations is a point of small height in the nullspace of the matrix of this linear system. Thus this principle can be viewed as a statement about points of small height in a given vector space. We write $H$ and $\H$ for appropriately selected height functions, which we will precisely define below. The following modern formulation of this result follows from a celebrated theorem of Bombieri and Vaaler, \cite{vaaler:siegel}.

\begin{thm}[\cite{vaaler:siegel}] \label{Siegel:gen} Let $K$ be a number field of degree $d$ and discriminant $\D_K$, and let $N \geq 1$ be an integer. Let $W$ be a non-zero subspace of $K^N$ of dimension $w \leq N$. There exists a non-zero point $\bx \in W$ such that
\begin{equation}
\label{Siegel:2}
H(\bx) \leq  \left\{ N |\D_K|^{1/d} \right\}^{1/2} \H(W)^{1/w}.
\end{equation}
\end{thm}

The exponent on $\H(W)$ in the upper bound of Theorem \ref{Siegel:gen} is best possible, however the constant is not. The best possible constant for Siegel's Lemma was recently obtained by Vaaler in \cite{vaaler:siegel_best}. The actual Bombieri - Vaaler theorem is more general: it produces a full basis of small height for $W$. Results of this sort were originally treated as important technical lemmas used in transcendental number theory and Diophantine approximations for the purpose of constructing a certain auxiliary polynomial (see \cite{vaaler:siegel} and \cite{bombieri:cohen} for more information). Nowadays they have evolved as important results in their own right.  
\bigskip

In this paper we consider a generalization of this problem. Let $K$ be a number field, and let $W$ be a subspace of $K^N$, $N \geq 2$. Let $V_1,...,V_M$ be subspaces of $K^N$ of dimension less than dimension of $W$. We want to prove the existence of a non-zero point of small height in $W \setminus \bigcup_{i=1}^M V_i$ providing an explicit upper bound on the height of such a point. More precisely, our main result reads as follows.

\begin{thm} \label{numberfield:main} Let $K$ be a number field of degree $d$ with discriminant $\D_K$. Let $N \geq 2$ be an integer, $l=\left[ \frac{N}{2} \right]$, and let $W$ be a subspace of $K^N$ of dimension $w$, $1 \leq w \leq N$. Let $1 \leq s < w$ be an integer, and let $V_1,...,V_M$ be nonzero subspaces of $K^N$ with $\max_{1 \leq i \leq M} \{ \dim_K(V_i) \} \leq s$. There exists a point $\bx \in W \setminus \bigcup_{i=1}^M V_i$ such that
\begin{equation}
\label{numberfield:bound}
H(\bx) \leq \C_{K,N}(w,s) \H(W)^d \left\{ \left( \sum_{i=1}^M \frac{1}{\H(V_i)^d} \right)^{\frac{1}{(w-s)d}} + M^{\frac{1}{(w-s)d+1}} \right\},
\end{equation}
where
\begin{equation}
\label{main:constant}
\C_{K,N}(w,s) = 2^{w(d+3)} |\D_K|^{\frac{w}{2}} \left( (wd)^w \binom{Nd}{ld}^{\frac{1}{2d}} \right)^{\frac{1}{w-s}}.
\end{equation}
\end{thm}

The dependence on $\H(W)$ in the upper bound of Theorem \ref{numberfield:main} is sharp at least in the case $K=\que$. Let $M=1$, and take $V_1$ to be a subspace of $W$ of dimension $w-1$ generated by the vectors corresponding to the first $w-1$ successive minima of $W$ with respect to an adelic unit cube. Then the smallest vector in $W \setminus V_1$ will be the one corresponding to the $w$-th successive minimum, and its height can be as large as a constant multiple of $\H(W)$: this is a consequence of the adelic version of Minkowski's successive minima theorem and the Bombieri - Vaaler version of Siegel's lemma (see \cite{vaaler:siegel}). 
\smallskip

We separately discuss a special case of our main result, which can be thought of as an inverse of Siegel's Lemma. Suppose that $W=K^N$, and let $L_1(\bX),...,L_M(\bX)$ be $M$ linear forms in $N$ variables with coefficients in $K$. Then we can prove the existence of a point $\bx$ in $K^N$ of relatively small height such that $L_i(\bx) \neq 0$ for every $i=1,...,M$ (i.e. $\bx$ is outside of the union of nullspaces of linear forms). This discussion generalizes some results presented in the companion paper \cite{me:classic} in the case $K=\que$ to any number field. In particular, Theorem \ref{numberfield:main} can be viewed as a generalization of Theorem 5.1 of \cite{me:classic}. Although we employ similar principles in the proof, the techniques and ideas of \cite{me:classic} are more elementary and combinatorial in nature. 
\smallskip

This paper is structured as follows. In section 2 we present a technical lemma on the problem of counting integer lattice points in a closed cube in $\real^N$. In section 3 we use this counting mechanism to prove Theorem \ref{numberfield:main}. In section 4 we discuss some interesting corollaries of this result.
\bigskip

We start with some notation. let $K$ be a number field of degree $d$ over $\que$, $O_K$ its ring of integers, $\D_K$ its discriminant, and $M(K)$ its set of places. For each place $v \in M(K)$ we write $K_v$ for the completion of $K$ at $v$ and let $d_v = [K_v:\que_v]$ be the local degree of $K$ at $v$, so that for each $u \in M(\que)$
\begin{equation}
\sum_{v \in M(K), v|u} d_v = d.
\end{equation}

\noindent
For each place $v \in M(K)$ we define the absolute value $\|\ \|_v$ to be the unique absolute value on $K_v$ that extends either the usual absolute value on $\real$ or $\cee$ if $v | \infty$, or the usual $p$-adic absolute value on $\que_p$ if $v|p$, where $p$ is a prime. We also define the second absolute value $|\ |_v$ for each place $v$ by $|a|_v = \|a\|_v^{d_v/d}$ for all $a \in K$. Then for each non-zero $a \in K$ the {\it product formula} reads
\begin{equation}
\label{product_formula}
\prod_{v \in M(K)} |a|_v = 1.
\end{equation} 

\noindent
For each finite place $v \in M(K)$, $v \nmid \infty$, we define the {\it local ring of $v$-adic integers} $O_v = \{ x \in K : |x|_v \leq 1 \}$, whose unique maximal ideal is $P_v =  \{ x \in K : |x|_v < 1 \}$. Then $O_K = \bigcap_{v \nmid \infty} O_v$.

\noindent
We extend absolute values to vectors by defining the local heights. For each $v \in M(K)$ define a local height $H_v$ on $K_v^N$ by
\begin{equation}
H_v(\bx) = \max_{1 \leq i \leq N} |x_i|_v,
\end{equation}
for each $\bx \in K_v^N$. Also, for each $v | \infty$ we define another local height
\begin{equation}
\H_v(\bx) = \left( \sum_{i=1}^N \|x_i\|_v^2 \right)^{d_v/2d} 
\end{equation}
Then we can define two slightly different global height functions on $K^N$:
\begin{equation}
H(\bx) = \prod_{v \in M(K)} H_v(\bx),\ \ \H(\bx) = \prod_{v \nmid \infty} H_v(\bx) \times \prod_{v | \infty} \H_v(\bx),
\end{equation}
for each $\bx \in K^N$. It is easy to see that
\begin{equation}
H(\bx) \leq \H(\bx) \leq \sqrt{N} H(\bx).
\end{equation}
All our inequalities will use height $H$ for vectors, however we use $\H$ to define the conventional Schmidt height on subspaces in the manner described below. This choice of heights coincides with \cite{vaaler:siegel}.
\smallskip

We extend both heights $H$ and $\H$ to polynomials by viewing them as height functions of the coefficient vector of a given polynomial. We also define a height function on subspaces of $K^N$. Let $V \subseteq K^N$ be a subspace of dimension $J$, $1 \leq J \leq N$. Choose a basis $\bx_1,...,\bx_J$ for $V$, and write $X = (\bx_1\ ...\ \bx_J)$ for the corresponding $N \times J$ basis matrix. Then 
$$V = \{ X \bt : \bt \in K^J \}.$$
On the other hand, there exists an $(N-J) \times N$ matrix $A$ with entries in $K$ such that 
$$V = \{ \bx \in K^N : A \bx = 0 \}.$$
Let $\I$ be the collection of all subsets $I$ of $\{1,...,N\}$ of cardinality $J$. For each $I \in \I$ let $I'$ be its complement, i.e. $I' = \{1,...,N\} \setminus I$, and let $\I' = \{ I' : I \in \I\}$. Then 
$$|\I| = \binom{N}{J} = \binom{N}{N-J} = |\I'|.$$
For each $I \in \I$, write $X_I$ for the $J \times J$ submatrix of $X$ consisting of all those rows of $X$ which are indexed by $I$, and $_{I'} A$ for the $(N-J) \times (N-J)$ submatrix of $A$ consisting of all those columns of $A$ which are indexed by $I'$. By the duality principle of Brill-Gordan \cite{gordan:1} (also see Theorem 1 on p. 294 of \cite{hodge:pedoe}), there exists a non-zero constant $\gamma \in K$ such that
\begin{equation}
\label{duality}
\det (X_I) = (-1)^{\varepsilon(I')} \gamma \det (_{I'} A),
\end{equation}
where $\varepsilon(I') = \sum_{i \in I'} i$. Define the vectors of {\it Grassmann coordinates} of $X$ and $A$ respectively to be 
$$Gr(X) = (\det (X_I))_{I \in \I} \in K^{|I|},\ \ Gr(A) = (\det (_{I'} A))_{I' \in \I'} \in K^{|I'|},$$
and so by (\ref{duality}) and (\ref{product_formula})
$$\H(Gr(X)) = \H(Gr(A)).$$
Define the height of $V$ denoted by $\H(V)$ to be this common value. This definition is legitimate, since it does not depend on the choice of the basis for $V$. In particular, notice that if 
$$L(X_1,...,X_N) = \sum_{i=1}^N q_i X_i \in K[X_1,...,X_N]$$
is a linear form with a non-zero coefficient vector $\bq \in K^N$, and $V = \{ \bx \in K^N : L(\bx) = 0 \}$ is an $(N-1)$-dimensional subspace of $K^N$, then
\begin{equation}
\label{1.4}
\H(V) = \H(L) = \H(\bq).
\end{equation}
\smallskip

The method of proof of Theorem \ref{numberfield:main} is the following. For a positive $R \geq 1$ we estimate cardinalities of sets
$$S_R(W) = \{ \bx \in W \cap O_K^N : \max_{v | \infty} H_v(\bx)^{d/d_v} \leq R \},$$
and $S_R(V_i) = S_R(W) \cap V_i$ for each $1 \leq i \leq M$. In other words, we count the number of points in sections of the adelic cube with ``sidelength'' $R$ by $W$ and by each $V_i$. Then we find $R$ large enough so that $|S_R(W)|$ is greater than $\sum_{i=1}^M |S_R(V_i)|$. A related estimate for the number of points of bounded height in a subspace of $K^N$ is provided by J. Thunder in \cite{thunder:number}. Thunder's estimate, however, is asymptotic with an implicit constant in the error term. This is not suitable for our purposes, since we need explicit upper and lower bounds. Our estimates are different from Thunder's also in the way that we are considering points inside of an adelic cube, which is a smaller set than the one considered in \cite{thunder:number}. We formulate our counting estimate precisely in Lemma \ref{adelic} at the end of section 3. We are now ready to proceed. Results of this paper also appear as a part of \cite{me:diss}. 
\bigskip

\section{Lattice points in cubes}

In this section we state some bounds on the number of points of a lattice in $\real^N$ inside of a closed cube. These will later be used to prove our main result. 
\smallskip

For the rest of this paper, let $R \geq 1$, and define
$$C^N_R = \{\bx \in \real^N : \max_{1 \leq i \leq N} |x_i| \leq R \},$$
to be a cube in $\real^N$ centered at the origin with sidelength $2R$. Given a lattice $\Lambda$ in $\real^N$ of rank $N$ and determinant $\Delta$, we want to estimate the quantity $|\Lambda \cap C_R^N|$. First suppose that $\rk(\Lambda)=N$. Then there exists an uppertriangular, nonsigular $N \times N$ matrix $A=(a_{mn})$ with positive real entries such that $\Lambda = \{ A \bxi : \bxi \in \zed^N \}$. Then by Corollary 3.3 of \cite{me:classic}, we have:
\begin{equation}
\label{box:lattice}
\prod_{m=1}^N \left[ \frac{2R}{a_{mm}} \right] \leq |\Lambda \cap (C_R^N+\bz)| \leq \prod_{m=1}^N \left( \left[ \frac{2R}{a_{mm}} \right] + 1 \right),
\end{equation}
for each point $\bz$ in $\real^N$. Notice that if $2R \geq \max_{1 \leq m \leq N} a_{mm}$, then the lower bound of (\ref{box:lattice}) is greater or equal than $\prod_{m=1}^N \left( \frac{2R}{a_{mm}} - 1 \right)$.
\smallskip

\noindent
If the matrix $A$ as above with fixed determinant $\Delta$ is such that all diagonal entries $a_{mm} \geq c$ for some positive constant $c$, then the right hand side of (\ref{box:lattice}) takes its maximum value and the left hand side takes its minimum value when $a_{mm} = c$ for $N-1$ distinct values of $m$. This leads to the following lemma.

\begin{lem} \label{2.4.1} Let $\Lambda$ be a lattice of full rank in $\real^N$ of determinant $\Delta$ such that there exists a positive constant $c$ and an uppertriangular basis matrix $A = (a_{mn})_{1 \leq m,n \leq N}$ of $\Lambda$ with diagonal entries $a_{mm} \geq c$ for all $1 \leq m \leq N$ (in particular, this is true with $c=1$ if $\Lambda \subseteq \zed^N$). Assume that $2R \geq \max \left\{ \frac{\Delta}{c^{N-1}}, c \right\}$. Then for each point $\bz$ in $\real^N$ we have
\begin{eqnarray}
\label{cube:lattice}
\left( \frac{2R c^{N-1}}{\Delta} - 1 \right)\left( \frac{2R}{c} - 1 \right)^{N-1} & \leq & |\Lambda \cap (C_R^N+\bz)| \nonumber \\
& \leq & \left( \frac{2R c^{N-1}}{\Delta} + 1 \right)\left( \frac{2R}{c} + 1 \right)^{N-1}.
\end{eqnarray}
\end{lem}

\noindent
Notice that the assumption on $R$ is not needed for the upper bound of (\ref{cube:lattice}). Moreover, this upper bound is sharp: consider the lattice $\Lambda = \Delta \zed \times \zed^{N-1}$ for a fixed $\Delta$. 
\bigskip

\section{Proof of Theorem \ref{numberfield:main}}

In fact, we prove a slightly sharper bound that reads as follows.

\begin{thm} \label{numberfield:main1} Let $K$ be a number field of degree $d$ with discriminant $\D_K$ and $r_2$ complex places. Let $N \geq 2$ be an integer, and let $W$ be a subspace of $K^N$ of dimension $w$, $1 \leq w \leq N$. Let $1 \leq s < w$ be an integer, and let $V_1,...,V_M$ be nonzero subspaces of $K^N$ of corresponding dimensions $l_1,...,l_M \geq 1$ with $\max_{1 \leq i \leq M} \{ l_i \} \leq s$. Define
\begin{equation}
\label{numberfield:bound1}
R_1 = \left( \left( \C^1_K(w) \H(W) \right)^{\frac{1}{w-s}} + 1 \right) \left\{ \left( \sum_{i=1}^M \frac{\C^2_{K,N}(l_i)}{\H(V_i)^d} \right)^{\frac{1}{(w-s)d}} + M^{\frac{1}{(w-s)d+1}} \right\},
\end{equation}
where
\begin{equation}
\label{numberfield:constant1}
\C^1_K(w) = 4^{\frac{w(2d-r_2)+1}{2d}} (wd)^{w} |\D_K|^{\frac{w}{2d}},\ \ \ \ \C^2_{K,N}(l_i) = \frac{2^{l_i r_2} \binom{Nd}{l_i d}^{1/2}}{|\D_K|^{l_i/2}},
\end{equation}
and
\begin{equation}
\label{numberfield:bound2}
R_2 = 2^{\frac{w(d-2r_2)}{2}} wd |\D_K|^{\frac{w}{2}} \H(W)^d.
\end{equation}
There exists a point $\bx \in W \setminus \bigcup_{i=1}^M V_i$ such that
$$H(\bx) \leq \max \{ R_1, R_2 \}.$$
\end{thm}

\proof
Let 
$$\sigma_1,...,\sigma_{r_1},\tau_1,...,\tau_{r_2},...,\tau_{2r_2}$$
be the embeddings of $K$ into $\cee$ with $\sigma_1,...,\sigma_{r_1}$ being real embeddings and $\tau_i,\tau_{r_2+i} = \bar{\tau}_i$ for each $1 \leq i \leq r_2$ being the pairs of complex conjugate embeddings. For each $\alpha \in K$ and each complex embedding $\tau_i$, write $\tau_{i1}(\alpha) = \Re(\tau_i(\alpha))$ and $\tau_{i2}(\alpha) = \Im(\tau_i(\alpha))$, where $\Re$ and $\Im$ stand respectively for real and imaginary parts of a complex number. We will view $\tau_i(\alpha)$ as a pair $(\tau_{i1}(\alpha), \tau_{i2}(\alpha)) \in \real^2$. Then $d=r_1+2r_2$, and for each $N \geq 1$ we define an embedding
$$\sigma^N=(\sigma_1^N,...,\sigma_{r_1}^N, \tau_1^N,...,\tau_{r_2}^N): K^N \longrightarrow K_{\infty}^N,$$
where
$$K_{\infty} = \prod_{v|\infty} K_v = \prod_{v|\infty} \real^{d_v} = \real^d,$$
since $\sum_{v|\infty} d_v = d$. Then $\sigma^N(O_K^N)$ can be viewed as a lattice of full rank in $\real^{Nd}$. 
\smallskip

For $R \geq 1$ let $C_R^{Nd}$ be the cube with sidelength $2R$ centered at the origin in $\real^{Nd}$, as above. Let $V$ be a subspace of $K^N$ of dimension $l$, $1 \leq l \leq N$. We want to estimate the number of lattice points in the slice of a cube by $\sigma^N(V)$. Let 
$$\Lambda(V) = \sigma^N \left( V \cap O_K^N \right),$$
then, by Theorem 2 of \cite{thunder:asymptotic}, $\Lambda(V)$ is a lattice in $\real^{Nd}$ of rank $ld$, and
\begin{equation}
\label{n2}
|\det(\Lambda(V))| = \left(\frac{|\D_K|^{1/2}}{2^{r_2}}\right)^l \H(V)^d.
\end{equation}
Notice that the exponent $d$ on $\H(V)$ appears because our height is absolute unlike the one in Theorem 2 of \cite{thunder:asymptotic}. Also, the constant $2^{-r_2}$ appears because we use a slightly different embedding into $\real^{Nd}$ than that in Theorem 2 of \cite{thunder:asymptotic} (see Lemma 2 on p. 115 of \cite{lang:ant}).
\smallskip

On the other hand, let $\bx_1,...,\bx_{ld}$ be a basis for $\Lambda(V)$ as a lattice in $\real^{Nd}$, and write $X = (\bx_1\ ...\ \bx_{ld}) = (x_{ij})$ for the $Nd \times ld$ basis matrix. Then each row of $X$ consists of blocks of all conjugates of $l$ algebraic integers from $O_K$.
If $I \subset \{ 1,...,Nd \}$ with $|I|=ld$, then write $X_I$ for the $ld \times ld$ submatrix of $X$ whose rows are rows of $X$ indexed by $I$. In other words, $X_I$ is the $I$-th Grassmann component matrix of $X$. Then each row of $X_I$ again consists of blocks of all conjugates of $l$ algebraic integers from $O_K$.
\smallskip

\noindent
Let $\{v_1,...,v_{r_1}\} \subset M(K)$ be the places corresponding to the real embeddings $\sigma_1,...,\sigma_{r_1}$, and let $\{u_1,...,u_{r_2}\} \subset M(K)$ be the places corresponding to the complex embeddings $\tau_1,...,\tau_{r_2}$. Let $\alpha \in O_K$, then $|\alpha|_v \leq 1$ for all $v \nmid \infty$, and so $|\alpha|_v \geq 1$ for at least one $v | \infty$, call this place $v_*$. If $v_*$ is real, say $v_*=v_j$ for some $1 \leq j \leq r_1$, then $|\sigma_j(\alpha)| \geq 1$. If $v_*$ is complex, say $v_*=u_j$ for some $1 \leq j \leq r_2$, then $\sqrt{\tau_{j1}(\alpha)^2 + \tau_{j2}(\alpha)^2} \geq 1$, hence $\max \{|\tau_{j1}(\alpha)|, |\tau_{j2}(\alpha)|\} \geq \frac{1}{\sqrt{2}}$. Therefore,
$$\max \{ |\sigma_1(\alpha)|,...,|\sigma_{r_1}(\alpha)|,|\tau_{11}(\alpha)|,|\tau_{12}(\alpha)|,...,|\tau_{r_2 1}(\alpha)|,|\tau_{r_2 2}(\alpha)| \} \geq \frac{1}{\sqrt{2}},$$
in other words the maximum of the Euclidean absolute values of all conjugates of an algebraic integer is at least $\frac{1}{\sqrt{2}}$. Therefore the maximum of the Euclidean absolute values of the entries of every row of $X_I$ is at least $\frac{1}{\sqrt{2}}$.
\smallskip

\noindent
By the Cauchy-Binet formula,
\begin{eqnarray}
\label{n3}
\max_{|I|=ld} |\det(X_I)|\  & \leq & \  |\det(\Lambda(V))| \nonumber \\
& = & \left( \sum_{|I|=ld} |\det(X_I)|^2 \right)^{1/2} \nonumber \\
& \leq & \binom{Nd}{ld}^{1/2} \max_{|I|=ld} |\det(X_I)|.
\end{eqnarray}
Let $J \subset \{ 1,...,Nd \}$ with $|J|=ld$ be such that $|\det(X_J)| = \max_{|I|=ld} |\det(X_I)|$, and let $\Omega(V)$ be the lattice of full rank in $\real^{ld}$ spanned over $\zed$ by the column vectors of $X_J$. By combining (\ref{n2}) and (\ref{n3}), we see that 
\begin{eqnarray}
\label{n4}
\binom{Nd}{ld}^{-1/2} \left(\frac{|\D_K|^{1/2}}{2^{r_2}}\right)^l \H(V)^d & = & \binom{Nd}{ld}^{-1/2} |\det(\Lambda(V))| \nonumber \\
& \leq & \det(\Omega(V)) = |\det(X_J)| \nonumber \\
& \leq & |\det(\Lambda(V))| = \left(\frac{|\D_K|^{1/2}}{2^{r_2}}\right)^l \H(V)^d.
\end{eqnarray}
For convenience, we denote $\det(\Omega(V))$ by $\Delta(V)$. By Corollary 1 on p. 13 of \cite{cass:geom}, we can select a basis for $\Omega(V)$ so that the basis matrix is upper triangular, all of its nonzero entries are positive, and the maximum entry of each row occurs on the diagonal. Each of these maximum values is at least $\frac{1}{\sqrt{2}}$, since each row still consists of blocks of all conjugates of $l$ algebraic integers from $O_K$. Therefore the lattice $\Omega(V)$ satisfies the conditions of Lemma \ref{2.4.1} with $c=\frac{1}{\sqrt{2}}$. Hence
\begin{equation}
\label{n5}
|\Omega(V) \cap C_R^{ld}| \leq \left( \frac{2^{\frac{3}{2}} R}{ 2^{\frac{ld}{2}} \Delta(V) }+1 \right) (2^{\frac{3}{2}} R+1)^{ld-1}.
\end{equation}
On the other hand, by Theorem 4.3 of \cite{me:classic} (in particular see equation (31) of \cite{me:classic}), we have
\begin{equation}
\label{n6.1}
|\Lambda(V) \cap C_R^{Nd}| \geq |\Omega(V) \cap C_{\frac{R}{ld}}^{ld}|.
\end{equation}
Assume that $R \geq 2^{\frac{ld}{2}} ld \Delta(V)$. Then combining (\ref{n6.1}) with the lower bound of Lemma \ref{2.4.1}, we obtain
\begin{eqnarray}
\label{n6}
|\Lambda(V) \cap C_R^{Nd}| & \geq & \left( \frac{2^{\frac{3}{2}} R}{ 2^{\frac{ld}{2}} ld \Delta(V) } - 1 \right) \left( \frac{2^{\frac{3}{2}} R}{ld} - 1 \right)^{ld-1} \nonumber \\
& \geq & \frac{1}{2^{\frac{ld}{2}} \Delta(V)} \left( \frac{R \left( 2^{\frac{3}{2}} - 1 \right)}{ld} \right)^{ld} \nonumber \\
& > & \frac{R^{ld}}{(ld)^{ld} \Delta(V)},
\end{eqnarray}
since $2^{\frac{3}{2}} - 1 > \frac{3}{2} > 2^{\frac{1}{2}}$.
\smallskip

For future use, we also need to define a projection $\varphi_V : \Lambda(V) \longrightarrow \Omega(V)$, given by our construction. Namely, if $X \bwy \in \Lambda(V)$ for some $\bwy \in \zed^{Nd}$, then $\varphi_V(X \bwy) = X_J \bwy_J$, where $\bwy_J \in \zed^{ld}$ is obtained from $\bwy$ by removing all the coordinates which are not indexed by $J$. It is quite easy to see that $\varphi_V$ is a $\zed$-module isomorphism.
\smallskip

Now let $W$ be a $w$-dimensional subspace of $K^N$, and let $V_1,...,V_M$ be $M$ proper subspaces of $W$ of respective dimensions $1 \leq l_1,...,l_M \leq s$. For $R \geq 1$, let
\begin{equation}
\label{srw}
S_R(W) = \{ \bx \in W \cap O_K^N : \max_{v | \infty} H_v(\bx)^{d/d_v} \leq R \},
\end{equation}
and for each $1 \leq i \leq M$, let $S_R(V_i) = S_R(W) \cap V_i$. Define a counting function
$$f_W(R) = |S_R(W)| - \left| \bigcup_{i=1}^M S_R(V_i) \right| \geq |S_R(W)| - \sum_{i=1}^M |S_R(V_i)|,$$
so that if $f_W(R)>0$ then there exists a point of height at most $R$ in $W \cap O_K^N$ outside of $\bigcup_{i=1}^M V_i$. Thus we want to find the minimal possible $R$ for which $f_W(R)>0$. 
\smallskip

\noindent
Notice that for each $\bx \in K^N$, 
$$\max_{v | \infty} H_v(\bx)^{d/d_v} = \max_{1 \leq j \leq N} \max \{ |\sigma_1(x_j)|,...,|\sigma_{r_1}(x_j)|,|\tau_1(x_j)|,...,|\tau_{r_2}(x_j)| \},$$
hence $\sigma^N(S_R(W)) = \sigma^N(W \cap O_K^N) \cap C_R^{Nd}$, and so $|S_R(W)| = |\sigma^N(S_R(W))| = |\Lambda(W) \cap C_R^{Nd}|$, since $\sigma^N$ is injective. Also, for each $1 \leq i \leq M$ the map $\varphi_{V_i} \circ \sigma^N$ is injective, and if for some $\bx \in S_R(V_i)$, $\bwy = \varphi_{V_i} \circ \sigma^N(\bx)$, then
$$R \geq \max_{v | \infty} H_v(\bx)^{d/d_v} \geq \max_{1 \leq j \leq l_i d} |y_j|,$$
therefore $\bwy \in \Omega(V_i) \cap C_R^{l_id}$. This means that for each $1 \leq i \leq M$, we have $|S_R(V_i)| \leq |\Omega(V_i) \cap C_R^{l_i d}|$. Hence we have proved that
$$f_W(R) \geq |\Lambda(W) \cap C_R^{Nd}| - \sum_{i=1}^M |\Omega(V_i) \cap C_R^{l_i d}|,$$
where the notation is as above. From here on assume that $R \geq 2^{\frac{wd}{2}} wd \Delta(W)$. Applying (\ref{n5}) and (\ref{n6}) we obtain
\begin{eqnarray}
\label{n7}
f_W(R) & \geq & \frac{R^{wd}}{(wd)^{wd} \Delta(W)} - \sum_{i=1}^M \left( \frac{R}{ 2^{\frac{l_i d-3}{2}} \Delta(V_i) }+1 \right) (2^{\frac{3}{2}} R+1)^{l_i d-1} \nonumber \\
                  & \geq & \frac{R^{wd}}{(wd)^{wd} \Delta(W)} - (2^{\frac{3}{2}} R+1)^{sd-1} \sum_{i=1}^M \left( \frac{R}{ 2^{\frac{d-3}{2}} \Delta(V_i) }+1 \right) \nonumber \\
                  & \geq & \frac{R^{wd}}{(wd)^{wd} \Delta(W)} - 4^{\left(s-\frac{1}{4}\right)d-\frac{1}{4}} \left( \sum_{i=1}^M \frac{1}{\Delta(V_i)} \right) R^{sd} - 4^{sd-1} M R^{sd-1} \nonumber \\
                  & \geq & \left( \frac{R^{sd-1}}{(wd)^{wd} \Delta(W)} \right) \times \nonumber \\
                  & \times & \left\{ R^{(w-s)d+1} - (4wd)^{wd} \Delta(W) \left( \sum_{i=1}^M \frac{1}{\Delta(V_i)} \right) R - (4wd)^{wd} \Delta(W) M \right\}.
\end{eqnarray}
Let $x = \sum_{i=1}^M \frac{1}{\Delta(V_i)}$, and let $\A_W = (4wd)^{wd} \Delta(W)$, and define
$$g_W(R) = R^{(w-s)d+1} - \A_W x R - \A_W M,$$
so that $f_W(R) \geq \frac{R^{sd-1}}{(wd)^{wd} \Delta(W)}  g_W(R)$. Hence we want to determine a value of $R$ for which $g_W(R)>0$. Let $\B_W$ be a positive number to be specified later. Then
\begin{eqnarray}
\label{n8-0}
\lefteqn{ g_W \left( \B_W \left( M^{\frac{1}{(w-s)d+1}}+x^{\frac{1}{(w-s)d}} \right) \right) } \nonumber \\
    & = & \B_W^{(w-s)d+1} \left( M^{\frac{1}{(w-s)d+1}}+x^{\frac{1}{(w-s)d}} \right)^{(w-s)d+1} \nonumber \\
    & & -\ \ \A_W \B_W \left( M^{\frac{1}{(w-s)d+1}}+x^{\frac{1}{(w-s)d}} \right) x - \A_W M \nonumber \\
    & \geq & (\B_W^{(w-s)d+1} - \A_W) M \nonumber \\
    & & +\ \ \B_W (\B_W^{(w-s)d} - \A_W) x^{1+\frac{1}{(w-s)d}} - \A_W \B_W M^{\frac{1}{(w-s)d+1}} \nonumber \\
    & \geq & (\B_W^{(w-s)d+1} - \A_W (\B_W+1)) M + \B_W (\B_W^{(w-s)d} - \A_W) x^{1+\frac{1}{(w-s)d}} \nonumber \\
    & > & 0,
\end{eqnarray}
for all $M$ and $x$ if $\B_W \geq 1$, and $\B_W^{(w-s)d} - 2 \A_W > 0$, hence we can choose 
\begin{eqnarray}
\label{n8}
\B_W & = & (2 \A_W)^{\frac{1}{(w-s)d}} + 1 = \left( 4^{wd+\frac{1}{2}} (wd)^{wd} \Delta(W) \right)^{\frac{1}{(w-s)d}} + 1 \nonumber \\
& \leq & \left( 4^{\frac{w(2d-r_2)+1}{2}} (wd)^{wd} |\D_K|^{\frac{w}{2}} \H(W)^d \right)^{\frac{1}{(w-s)d}} + 1,
\end{eqnarray}
where the last inequality follows by (\ref{n4}). Therefore, $f_W(R)>0$ if $R$ is such that
\begin{eqnarray}
\label{n9}
R & \geq & \left\{ \left( 4^{\frac{w(2d-r_2)+1}{2}} (wd)^{wd} |\D_K|^{\frac{w}{2}} \H(W)^d \right)^{\frac{1}{(w-s)d}} + 1 \right\} \times \nonumber \\
& \times & \left\{ \left( \sum_{i=1}^M \frac{1}{\Delta(V_i)} \right)^{\frac{1}{(w-s)d}} + M^{\frac{1}{(w-s)d+1}} \right\}.
\end{eqnarray}
Estimating the latter from above using (\ref{n4}), we infer that $f_W(R)>0$ if
\begin{eqnarray}
\label{n9.1}
 R & \geq & \left\{ \left( 4^{\frac{w(2d-r_2)+1}{2}} (wd)^{wd} |\D_K|^{\frac{w}{2}} \H(W)^d \right)^{\frac{1}{(w-s)d}} + 1 \right\} \times \nonumber \\
& \times & \left\{ \left( \sum_{i=1}^M \frac{2^{l_i r_2} \binom{Nd}{l_i d}^{1/2}}{|\D_K|^{l_i/2} \H(V_i)^d} \right)^{\frac{1}{(w-s)d}} + M^{\frac{1}{(w-s)d+1}} \right\}.
\end{eqnarray}
By our original assumption $R$ must also be greater or equal than $2^{\frac{wd}{2}} wd \Delta(W)$. To accomplish this, by (\ref{n4}) we can take
\begin{equation}
\label{n9.2}
R \geq 2^{\frac{w(d-2r_2)}{2}} wd |\D_K|^{\frac{w}{2}} \H(W)^d.
\end{equation}
Combining (\ref{n9.1}) with (\ref{n9.2}) completes the proof.
\endproof
\bigskip

Notice that the main part of this argument can be treated as a separate result on the number of points of a subspace of $K^N$ in the adelic cube. Write $K_{\aaa}$ for the ring of the adeles of $K$. Define the $N$-dimensional adelic cube with ``sidelength'' $R$ to be
\begin{equation}
C_{\aaa}^N (R) = \prod_{v \nmid \infty} O^N_v \times \prod_{v | \infty} \{ \bx \in K^N_v : H_v(\bx)^{d/d_v} \leq R \},
\end{equation}
for $R \geq 1$. This is a basic example of a compact convex symmetric set in the adelic geometry of numbers (see \cite{vaaler:siegel} for details). $K^N$ can be viewed as a lattice in $K_{\aaa}^N$ under the standard diagonal embedding. For a subspace $W$ of $K^N$ we also write $W$ for its image under this embedding. Clearly $C_{\aaa}^N (R) \cap W$ is a finite set. In fact, it is precisely the set $S_R(W)$ as defined by (\ref{srw}). The following lemma follows from the argument in the proof of Theorem \ref{numberfield:main1} above.

\begin{lem} \label{adelic} Let $W \subseteq K^N$ be a $w$-dimensional subspace, $1 \leq w \leq N$, and let $R \geq 1$. Then
\begin{eqnarray}
\label{adelic:bound}
\lefteqn{ \left( \frac{2^{\frac{w(2r_2-d)+3}{2}} R}{ wd |\D_K|^{\frac{w}{2}} \H(W)^d } - 1 \right) \left( \frac{2^{\frac{3}{2}} R}{wd} - 1 \right)^{wd-1} \leq \ \ \ |C_{\aaa}^N (R) \cap W| } \nonumber \\
& & \ \ \ \ \ \ \ \ \ \ \ \ \ \ \ \ \ \ \ \ \ \ \ \leq \left( \frac{\binom{Nd}{wd}^{\frac{1}{2}} 2^{\frac{w(2r_2 - d) + 3}{2}} R}{|\D_K|^{\frac{w}{2}} \H(W)^d} + 1 \right) ( 2^{\frac{3}{2}} R + 1 )^{wd-1}.
\end{eqnarray}
\end{lem}
\noindent
Lemma \ref{adelic} presents the counting principle that is our main tool.
\bigskip

\section{Corollaries}

Notice that in case $K = \que$ and $s=w-1$ the bound of Theorem \ref{numberfield:main} becomes
\begin{equation}
\label{n10}
(16 w)^w \binom{N}{l}^{1/2} \H(W) \left\{ \sum_{i=1}^M \frac{1}{\H(V_i)} + \sqrt{M} \right\},
\end{equation}
which is essentially (up to a constant) the bound of Theorem 5.1 in \cite{me:classic}. 
\smallskip

Here is another interesting observation that generalizes some ideas of \cite{me:classic}. Suppose that $W=K^N$ and $V_1,...,V_M$ is a collection of nullspaces of linear forms $L_1,...,L_M$ in $N$ variables with coefficients in $K$ (i.e. $w=N$ and $l_i=s=N-1$ for each $1 \leq i \leq M$). Let
$$F(X_1,...,X_N) = \prod_{i=1}^M L_i(X_1,...,X_N).$$
Then $F$ is a homogeneous polynomial of degree $M$ in $N$ variables with coefficients in $K$. Hence Theorem \ref{numberfield:main1} produces a point $\bx \in K^N$ of small height at which $F$ does not vanish. In fact, a simple explicit bound on $H(\bx)$ that depends only on $K$, $N$, and $M$ follows from Theorem \ref{numberfield:main1} in this case:
\begin{equation}
\label{n11}
H(\bx) \leq 2^{N(d+3)+1} \left( Nd |\D_K| \right)^{\frac{N}{2}} \binom{Nd}{Nd-d}^{\frac{1}{2d}} M^{1/d}.
\end{equation}
Notice that this is a certain inverse of Siegel's Lemma: we produce a point of small height outside of a collection of subspaces. This can also be viewed as an effective instance of the following more general non-effective simple lemma.

\begin{lem} \label{2.2.3} Let $K$ be a number field of degree $d$, and let $F$ be a polynomial in $N \geq 2$ variables of degree $M \geq 1$ with coefficients in $K$. There exists a constant $\C_K(N)$ and $\bx \in O_K^N$ such that $F(\bx) \neq 0$, and
\begin{equation}
\label{bound:Md}
H(\bx) \leq \C_K(N) M^{1/d}.
\end{equation}
\end{lem}
 
\proof
Let
$$S_M(K) = \left\{ x \in K :|x|_v \leq 1\ \forall\ v \nmid \infty,\ \ |x|^{d/d_v}_v \leq \C(K) M^{1/d}\ \forall\ v | \infty \right\},$$
where $\C(K)$ is a positive field constant to be specified later. By \cite{lang:ant} (Theorem 0, p. 102) there exist constants $\A(K)$ and $\B(K)$ such that
\begin{equation}
\label{lang:p102}
\A(K) \C(K)^d M \leq |S_M(K)| \leq \B(K) \C(K)^d M.
\end{equation}
Let
\begin{equation}
\C(K) = \left( \frac{2}{\A(K)} \right)^{1/d},
\end{equation}
so that $|S_M(K)| \geq 2M \geq M+1$. It is a well-known fact (see for instance Lemma 1 on p. 261 of \cite{cass:geom}, also Lemma 2.1 of \cite{me:classic}) that a non-zero polynomial of degree $M$ in $N$ variables cannot vanish on the whole set $S^N$ if $S$ is a set of cardinality larger than $M$. Hence there must exist $\bx \in S_M(K)^N$ such that $F(\bx) \neq 0$, and so
\begin{equation}
H(\bx) \leq \prod_{v | \infty} \left( \C(K) M^{1/d} \right)^{d_v/d} = \C(K) M^{1/d}.
\end{equation}
This completes the proof.
\endproof

Notice that the upper bound in (\ref{bound:Md}) has the correct order of magnitude in the following sense. It is conceptual for the cardinality of the set $S_M(K)$ in the proof of Lemma \ref{2.2.3} to be at least $M+1$, since  there are polynomials of degree $M$ that vanish on a set $S^N$ if $|S| \leq M$: let $S = \{\alpha_1,...,\alpha_M\} \subset \zed$, and let
$$F(X_1,...,X_N) = \sum_{i=1}^N \prod_{j=1}^M (X_i - \alpha_j).$$
\bigskip

Another interesting immediate corollary of Theorem \ref{numberfield:main} in the case $M=1$ is the following subspace extension lemma.

\begin{cor} \label{extend:subspace} Let $K$ be a number field as in Theorem \ref{numberfield:main1}. Let $N \geq 2$ be an integer, and let $W$ be a subspace of $K^N$ of dimension $w$, $1 < w \leq N$. Let $V \subseteq W$ be a proper subspace of $W$ of dimension $(w-1) \geq 1$. There exists a point $\bx \in O_K^N$ such that $W = \spn_K \{V, \bx\}$, and
\begin{equation}
\label{extend:bound}
H(\bx) \leq \C_{K,N}(w,w-1) \H(W)^d \left( 1 + \frac{1}{\H(V)} \right),
\end{equation}
where the constant $\C_{K,N}(w,w-1)$ is as in (\ref{main:constant}).
\end{cor}
\smallskip

{\bf Aknowledgements.} I want to thank Professor Jeffrey D. Vaaler for his valuable advice and numerous useful conversations on the subject of this paper. I would also like to thank Professor Preda Mihailescu and the referee for their helpful comments.

\nocite{*}
\bibliographystyle{plain}  
\bibliography{fukshansky}        

\begin{thebibliography}{10}

\bibitem{bombieri:cohen}
E.~Bombieri and P.~B. Cohen.
\newblock Siegel's lemma, {P}ade approximations and {J}acobians.
\newblock {\em Ann. Scuola Norm. Sup. Pisa Cl. Sci. (4)}, 25(1-2):155--178,
  1998.

\bibitem{vaaler:siegel}
E.~Bombieri and J.~D. Vaaler.
\newblock On {S}iegel's lemma.
\newblock {\em Invent. Math.}, 73(1):11--32, 1983.

\bibitem{cass:geom}
J.~W.~S. Cassels.
\newblock {\em An Introduction to the Geometry of Numbers}.
\newblock Springer-Verlag, 1959.

\bibitem{me:classic}
L.~Fukshansky.
\newblock Integral points of small height outside of a hypersurface.
\newblock {\em to appear in Monatsh. Math.}

\bibitem{me:diss}
L.~Fukshansky.
\newblock {\em Algebraic points of small height with additional arithmetic
  conditions}.
\newblock PhD thesis, {U}niversity of {T}exas at {A}ustin, 2004.

\bibitem{gordan:1}
P.~Gordan.
\newblock Uber den grossten gemeinsamen factor.
\newblock {\em Math. Ann.}, 7:443--448, 1873.

\bibitem{hodge:pedoe}
W.~V.~D. Hodge and D.~Pedoe.
\newblock {\em Methods of Algebraic Geometry, Volume 1}.
\newblock Cambridge Univ. Press, 1947.

\bibitem{lang:ant}
S.~Lang.
\newblock {\em Algebraic Number Theory}.
\newblock Addison-Wesley, 1970.

\bibitem{siegel}
C.~L. Siegel.
\newblock Uber einige {A}nwendungen diophantischer {A}pproximationen.
\newblock {\em Abh. der Preuss. Akad. der Wissenschaften Phys.-math Kl.}, Nr.
  1:209--266, 1929.

\bibitem{thue}
A.~Thue.
\newblock Uber {A}nnaherungswerte algebraischer {Z}ahlen.
\newblock {\em J. Reine Angew. Math.}, 135:284--305, 1909.

\bibitem{thunder:asymptotic}
J.~L. Thunder.
\newblock An asymptotic estimate for heights of algebraic subspaces.
\newblock {\em Trans. Amer. Math. Soc.}, 331:395--424, 1992.

\bibitem{thunder:number}
J.~L. Thunder.
\newblock The number of solutions of bounded height to a system of linear
  equations.
\newblock {\em J. Number Theory}, 43:228--250, 1993.

\bibitem{vaaler:siegel_best}
J.~D. Vaaler.
\newblock The best constant in {S}iegel's lemma.
\newblock {\em Monatsh. Math.}, 140(1):71--89, 2003.

\end{thebibliography}

\end{document}